\documentclass[11pt]{article}
\usepackage{times}
\usepackage{indentfirst}
\usepackage{amsthm,amsmath}
\usepackage{calrsfs}

\newtheorem{Theorem}{Theorem}
\newtheorem{Lemma}{Lemma}

\newtheorem{Corollary}{Corollary}

\newcommand\mysection[1]{
              \refstepcounter{section}
              \section*{\normalsize\bf\thesection.~#1}
                       }

\begin{document}

\section*{\centering\Large\rm\bf
          ASYMPTOTICS OF RANDOMLY STOPPED SUMS\\
          IN THE PRESENCE OF HEAVY TAILS}
\section*{\centering\large\sc
By Denis Denisov,\footnote{Supported by EPSRC grant No.~EP/E033717/1.}
Sergey Foss,$^{1,}$\footnote{Supported
by the Royal Society International Joint Project Grant 2005/R2 UP.}
and Dmitry Korshunov$^{2,}$\footnote{Supported by
the Edinburgh Mathematical Society.}}
\section*{\centering\normalsize\it
Heriot-Watt University and Sobolev Institute of Mathematics}

\begin{abstract}
We study conditions under which
$$
{\bf P}\{S_\tau>x\}\sim {\bf P}\{M_\tau>x\}
\sim {\bf E}\tau{\bf P}\{\xi_1>x\}\ \mbox{ as }x\to\infty,
$$
where $S_\tau$ is a sum $\xi_1+\ldots+\xi_\tau$
of random size $\tau$ and $M_\tau$ is a maximum
of partial sums $M_\tau=\max_{n\le\tau}S_n$.
Here $\xi_n$, $n=1$, $2$, \ldots, are independent identically
distributed random variables whose common distribution
is assumed to be subexponential.
We consider mostly the case where $\tau$ is independent
of the summands; also, in a particular situation,
we deal with a stopping time.

Also we consider the case where ${\bf E}\xi>0$ and where
the tail of $\tau$ is comparable with or heavier than
that of $\xi$, and obtain the asymptotics
$$
{\bf P}\{S_\tau>x\}\sim
{\bf E}\tau {\bf P}\{\xi_1>x\}+{\bf P}\{\tau>x/{\bf E}\xi\}
\ \mbox{ as }x\to\infty.
$$
This case is of a primary interest in the branching processes.

In addition, we obtain new uniform (in all $x$ and $n$)
upper bounds for the ratio
${\bf P}\{S_n>x\}/{\bf P}\{\xi_1>x\}$ which
substantially improve Kesten's bound in the
subclass ${\mathcal S}^*$ of subexponential distributions.

{\it AMS 2000 subject classifications}: Primary 60E05;
secondary 60F10, 60G70

{\it Key words and phrases:}
Random sums of random variables;
Upper bound;
Convolution equivalence;
Heavy-tailed distributions;
Subexponential distributions
\end{abstract}

\hspace{20mm}

\mysection{Introduction}\label{intro}

Let $\xi$, $\xi_1$, $\xi_2$, \ldots~be independent
identically distributed random variables with
a finite mean.
We assume that their common distribution
$F$ is right-unbounded, that is,
$\overline F(x) \equiv{\bf P}\{\xi>x\}>0$ for all $x$.
Moreover, we assume that $F$ has a {\it heavy} (right) tail.
Recall that a random variable $\eta$ has a {\it heavy-tailed}
distribution
if ${\bf E} e^{\varepsilon\eta} = \infty$ for all
$\varepsilon>0$, and {\it light-tailed} otherwise.

Let $S_0=0$ and $S_n=\xi_1+\ldots+\xi_n$,
$n=1$, $2$, \ldots, and let
$M_n= \max_{0\le i \le n}S_{i}$ be the partial maxima.
Denote by $F^{*n}$ the distribution of $S_n$.

Let $\tau$ be a counting random variable with a finite mean.
In this paper, we study the asymptotics
for the tail probabilities
${\bf P}\{S_\tau>x\}$ and ${\bf P}\{M_\tau>x\}$
as $x\to\infty$.

It is known that, for {\it any} distribution $F$ on ${\bf R}^+$
and for {\it any} counting random variable $\tau$ which
is independent of the sequence $\{\xi_n\}$,
$$
\liminf_{x\to\infty} \frac{{\bf P}\{S_\tau>x\}}{\overline F(x)}
\ge {\bf E}\tau,
$$
see, e.g. [\ref{Rudin}, \ref{DFK2}].
It was proved in the series of papers
[\ref{FK}, \ref{DFK1}, \ref{DFK2}] that
if $F$ is a heavy-tailed distribution on ${\bf R}^+$
with finite mean and if ${\bf P}\{c\tau>x\}=o(\overline F(x))$
as $x\to\infty$, for some $c>{\bf E}\xi$, then
\begin{eqnarray}\label{lim.inf}
\liminf_{x\to\infty} \frac{{\bf P}\{S_\tau>x\}}{\overline F(x)}
= {\bf E}\tau.
\end{eqnarray}
This gives us the idea what asymptotic behaviour of
${\bf P}\{S_\tau>x\}$ should be expected, at least
if the tail of $\tau$ is lighter than that of $\xi$.
In particular, by considering the case
$\tau=2$, we conclude that if $F$ is a heavy-tailed
distribution on ${\bf R}^+$ and if
${\bf P}\{S_2>x\}\sim c\overline F(x)$
as $x\to\infty$, for some $c$,
then $c=2$ with necessity (see [\ref{FK}]).
By the latter observation, we restrict our attention
to subexponential distributions only.

A distribution $F$ on ${\bf R}^+$ with unbounded
support is called {\it subexponential},
$F\in{\mathcal S}$, if
$\overline{F*F}(x) \sim 2\overline F(x)$
as $x\to\infty$.
A distribution $F$ on ${\bf R}$ is called
subexponential if its conditional distribution
on ${\bf R}^+$ is subexponential.
It is well known that any subexponential distribution
is heavy-tailed and, even more, is long-tailed.
A distribution $F$ with right-unbounded
support is called {\it long-tailed} if
$\overline F(x+y)\sim\overline F(x)$ as $x\to\infty$,
for any fixed $y$.

The key result in the theory of subexponential
distributions is: if $F$ is subexponential and
if $\tau$ does not depend on the summands and
is light-tailed, then
\begin{eqnarray}\label{tau.equi.ind}
{\bf P}\{S_\tau>x\} \sim
{\bf E}\tau\overline F(x)
\ \mbox{ as }x\to\infty.
\end{eqnarray}
A converse result also holds: if, for a distribution
$F$ on ${\bf R}^+$ and for an independent
counting random variable $\tau\ge 2$,
${\bf P}\{S_\tau>x\} \sim {\bf E}\tau\overline F(x)$
as $x\to\infty$, then $F$ is subexponential
(see, e.g. [\ref{EGV}]).

The intuition behind relation (\ref{tau.equi.ind})
is the {\it principle of one big jump}:
in the case of heavy tails, for $x$ large,
the most probable way leading to the event
$\{S_n>x\}$ is that one of $n$ summands
$\xi_1$, \ldots, $\xi_n$ is large while all others
are relatively small. Asymptotically this gives
the probability $n\overline F(x)$, and conditioning
on $\tau$ yields to the multiplier ${\bf E}\tau$.
The keystone of the proof is Kesten's bound:
for any subexponential distribution $F$ and for any
$\varepsilon>0$, there exists $K=K(F,\varepsilon)$
such that the inequality
$$
\overline{F^{*n}}(x) \le K(1+\varepsilon)^n\overline F(x)
$$
holds for all $x$ and $n$; see, e.g. [\ref{AN}, Section IV.4].
Clearly this estimate does not help
to prove (\ref{tau.equi.ind})
if the distribution of $\tau$ is heavy-tailed.
So the question of the basic importance is:
If we fix a subexponential distribution $F$,
then what are the weakest natural conditions on
$\tau$ which still guarantee relation
(\ref{tau.equi.ind}) to hold?
Intuitively, the light-tailedness assumption
seems to be very strong. The study of this problem
is one of the main topics of the present paper.

In order to state our first result, we need to introduce
the notion of ${\mathcal S}^*$-distribution.
A distribution $F$ on ${\bf R}$
with a finite mean belongs to the class ${\mathcal S}^*$ if
\begin{eqnarray*}
\int_0^x \overline F(x-y)\overline F(y)dy
\sim 2a\overline F(x)\quad\mbox{ as }x\to\infty,
\end{eqnarray*}
where $a=2\int_0^\infty\overline F(y)dy$.
It is known (see Kl\"uppelberg [\ref{Kl}])
that any distribution from the class
${\mathcal S}^*$ is subexponential.
Though these two classes, ${\mathcal S}^*$ and
${\mathcal S}$, are considered as rather similar,
there exist subexponential distributions
which are not in ${\mathcal S}^*$,
see, e.g. [\ref{DFK2004}] and the discussion
in Section \ref{bounds} below.
Classical examples of distributions from the class
${\mathcal S}^*$ are Pareto, log-normal, and Weibull with
parameter $\beta \in (0,1)$.

\begin{Theorem}\label{thm1}
Assume that a counting random variable $\tau$ does not
depend on $\{\xi_n\}$. Let $F\in{\mathcal S}^*$.

{\rm(i)} If ${\bf E}\xi<0$ then
\begin{eqnarray}\label{eq_lim}
{\bf P}\{S_\tau>x\} \sim {\bf P}\{M_\tau>x\}
\sim {\bf E}\tau\overline F(x)
\ \mbox{ as } x\to\infty.
\end{eqnarray}

{\rm(ii)} If ${\bf E}\xi\ge0$ and if
there exists $c>{\bf E}\xi$ such that
\begin{eqnarray}\label{eq1}
{\bf P}\{c\tau>x\} = o(\overline F(x))
\ \mbox{ as }x\to\infty,
\end{eqnarray}
then asymptotics {\rm(\ref{eq_lim})} again hold.
\end{Theorem}

The latter theorem shows that if we restrict our
attention from the class of all heavy-tailed
distributions to the class ${\mathcal S^*}$,
we obtain equivalence (\ref{eq_lim})
which is stronger than assertion (\ref{lim.inf})
for the `$\liminf$'. Definitely we should assume the
subexponentiality of $F$ in order to get (\ref{eq_lim}).
At the end of Section \ref{pth1} we construct
an example demonstrating that the stronger condition
$F\in {\mathcal S}^*$ is essential for the statement to hold in
the whole generality and cannot be replaced by condition
$F\in {\mathcal S}$.

The proof of Theorem \ref{thm1} is carried out
in Section \ref{pth1}. Statement (i)
can be found in [\ref{FZ_2003}]; in Section
\ref{pth1} we give an alternative proof of (i).
Note that these two cases,
negative and positive mean of $\xi$,
are substantially different in their nature.

Condition (\ref{eq1}) seems to be essential,
since, for any $c<{\bf E}\xi$,
\begin{eqnarray*}
{\bf P}\{S_\tau>x\}
&=& {\bf P}\{S_\tau>x,c\tau\le x\}
+{\bf P}\{S_\tau>x,c\tau>x\}\\
&\ge& ({\bf E}\tau+o(1))\overline F(x)
+(1+o(1)){\bf P}\{c\tau>x\}
\end{eqnarray*}
as $x\to\infty$, due to the convergence
${\bf P}\{S_\tau>x|c\tau>x\}\to 1$,
by the Law of Large Numbers.
In particular, for $\tau$ with a regularly
varying tail distribution, condition
(\ref{eq1}) is necessary
for asymptotic relation (\ref{eq_lim}) to hold.
Further discussion on condition (\ref{eq1})
can be found in Section \ref{pth1}.

Stam in [\ref{Stam}, Theorem 5.1] and
A. Borovkov and K. Borovkov in [\ref{BB}, Section 7.1]
obtained asymptotics (\ref{eq_lim}) under
condition (\ref{eq1}) for regularly varying $F$.
Some results by Stam [\ref{Stam}] have been proved
again by Fa\"y {\it et al.} in [\ref{FGMS}].
The case where $F$ is a dominated varying
distribution was studied by Ng {\it et al}. [\ref{Ng_etal}]
and by Daley {\it et al}. [\ref{DOV}].
A subclass of so-called semi-exponential $F$
was considered in [\ref{BB}, Section 7.2].
In [\ref{FZ_2003}, Corollary 2],
asymptotics (\ref{eq_lim}) were obtained
in the case ${\bf E}\xi\ge0$ under the extra
assumption ${\bf P}\{\tau>h(x)\}=o(\overline F(x))$
for some function $h(x)\to\infty$ such that
$\overline F(x\pm h(x))\sim\overline F(x)$.

In Section \ref{bounds}, we derive new simple
uniform upper bounds for the ratio
$\overline{F^{*n}}(x)/\overline F(x)$
which generalise Kesten's bound for
${\mathcal S}^*$-distributions.
We prove the following

\begin{Theorem}\label{th.asymp}
Assume that $F\in{\mathcal S}^*$.
If ${\bf E}\xi<0$,
then there exists a constant $K$ such that
\begin{eqnarray*}
\frac{\overline{F^{*n}}(x)}{\overline F(x)}
\le Kn \quad\mbox{ for all }n\mbox{ and }x.
\end{eqnarray*}
If ${\bf E}\xi\in[0,\infty)$,
then, for any $c>{\bf E}\xi$, there exists $K$ such that
\begin{eqnarray*}
\frac{\overline{F^{*n}}(x)}{\overline F(x)}
\le \frac{K}{\overline F(cn)}
\quad\mbox{ for all }n\mbox{ and }x.
\end{eqnarray*}
\end{Theorem}

The latter estimates are also of their own interest.
They substantially improve similar bounds
in Shneer [\ref{Shneer}, Theorems 1 and 2]
(see also Daley {\it et al}. [\ref{DOV}, Theorem 3]).
In Theorem \ref{not.n.tails}, Section \ref{bounds},
we show that the condition $F\in{\mathcal S}^*$
is essential for the statement of Theorem
\ref{th.asymp} to hold; more precisely, we construct
a distribution $F\in{\mathcal S}\setminus{\mathcal S}^*$
with negative mean such that
$\sup_{n,x}\frac{\overline{F^{*n}}(x)}
{n\overline F(x)}=\infty$.

A closely related topic is the asymptotics of the type
${\bf P}\{S_n>x\}\sim n\overline F(x)$ as $n$, $x\to\infty$
which have been extensively studied starting from 60s.
The first works are remarkable papers of
S. Nagaev [\ref{SN1960}, \ref{SN1962}],
Linnik [\ref{Linnik}] (in this paper, in a special case,
the asymptotics are stated,
but the key relation (10.10) on p. 303
is not supported by a proof),
and later on of A. Nagaev [\ref{AN1}, \ref{AN2}]
where in particular the regularly varying
distributions were considered.
Namely, if $F$ is regularly varying with
the parameter $\alpha>2$ and ${\bf E}\xi_1=0$,
${\bf E}\xi_1^2=1$, then under mild technical conditions
(see [\ref{AN1}], [\ref{SN}, Theorem 1.9], or
[\ref{Roz1990}, Theorem 6]) the following asymptotics hold
\begin{eqnarray*}
{\bf P}\{S_n>x\} \sim
\overline\Phi(x/\sqrt n)+n\overline F(x)
\quad\mbox{as }x\to\infty\mbox{ uniformly in }n\le x^2;
\end{eqnarray*}
here $\overline\Phi$ is the tail function
of the standard normal law. Further, it follows that,
if $x{\le}\sqrt{(\alpha{-}2{-}\varepsilon)n\ln n}$, then
the asymptotics follow the Cental Limit Theorem,
while if $x>\sqrt{(\alpha-2+\varepsilon)n\ln n}$,
then the probability of a single big jump dominates.
For Weibull-type distributions the situation
is more complicated, see, e.g., A. Nagaev [\ref{AN2}],
S. Nagaev [\ref{SN1973}], Rosovskii [\ref{Roz1993}, \ref{Roz1997}].
Detailed overviews of results
in the theory of large deviations for random walks
with subexponential increments are given
in S. Nagaev [\ref{SN}] and in
Mikosch and A. Nagaev [\ref{MAV}].
There is still an ongoing research in this area,
see recent works by A. Borovkov and K. Borovkov [\ref{BB}],
A. Borovkov and Mogulskii [\ref{BM}],
Denisov {\it et al.} [\ref{DDS}]
and references therein.
In Section \ref{nFx} of this paper,
for an {\it arbitrary} distribution
$F\in{\mathcal S}^*$, we find
a range for $n=n(x)$ where the asymptotics
${\bf P}\{S_n>x\}\sim n\overline F(x)$ hold.
The corresponding proof is surprisingly short.

In Section \ref{tailcom}, we study the case where the
tail distributions of $\tau$ and $\xi$ are asymptotically
comparable and, for a subclass of subexponential distributions,
we obtain the asymptotics for ${\bf P}\{S_\tau>x\}$
which differ from (\ref{eq_lim}), see Theorem \ref{thm2.D}.
This generalises results by A. Borovkov and K. Borovkov
[\ref{BB}] and by Stam [{\ref{Stam}],
see Section \ref{tailcom} for further comments.
As a corollary, in Section \ref{branch},
we obtain new tail asymptotics for
Galton--Watson branching processes.

In Section \ref{taudep}, we study the case
where $\tau$ may depend on $\{ \xi_n\}$,
in particular, where $\tau$ is a stopping time.
First, we prove Theorem \ref{thm.stopping}
where we obtain equivalence (\ref{eq_lim})
for bounded $\tau$.
In the proof, we adapt the approach developed in
[\ref{Greenwood}] and generalise Greenwood's
result onto the whole class of subexponential
distributions.
Then we consider an unbounded $\tau$ and prove
Theorem \ref{thm.stopping.} which states that
equivalence (\ref{eq_lim}) holds under
a stronger assumption than (\ref{eq1})
(see condition (\ref{hx.cond}) below).
Theorem \ref{thm.stopping.} geleralises earlier
results by Greenwood and Monroe [\ref{GreenwoodMonroe}]
and by A. Borovkov and Utev [\ref{BU}],
see Corollary \ref{cor.stopping.reg} and comments after it.
Concerning the asymptotics for the maximum, it was shown in
[\ref{FZ_2003}] (see also [\ref{FPZ}]) that the equivalence
${\bf P}\{M_\tau>x\}\sim{\bf E}\tau\overline F(x)$
holds without any further assumptions on the tail
distribution of $\tau$ if ${\bf E}\xi<0$
and under condition (\ref{hx.cond}) otherwise.

\mysection{Uniform upper bounds for tails;
proof of Theorem \ref{th.asymp}}\label{bounds}

In this Section,
for the ratios $\overline{F^{*n}}(x)/\overline F(x)$,
we derive more precise upper bounds
than Kesten's bound, which are again uniform in $x$.
We consider two cases ${\bf E}\xi<0$
and ${\bf E}\xi\ge 0$ separately.
We need the following result:

\begin{Theorem}[{\rm[\ref{K2002}]
and [\ref{DFK2004}, Corollary 4]}]\label{prop.maxima}
Assume that $F\in{\mathcal S}^*$ and ${\bf E}\xi<0$.
Then, as $x\to\infty$ and uniformly in $n\ge1$,
\begin{eqnarray*}
{\bf P}\{M_n>x\} \sim \frac{1}{|{\bf E}\xi|}
\int_x^{x+n|{\bf E}\xi|} \overline F(y)dy.
\end{eqnarray*}
\end{Theorem}

\proof[Proof of Theorem \ref{th.asymp}]
First we consider the case (i) of negative mean.
Taking into account the inequality
$S_n \le M_n$, Theorem \ref{prop.maxima},
and the inequality
\begin{eqnarray}\label{ine.n}
\frac{1}{|{\bf E}\xi|}
\int_x^{x+n|{\bf E}\xi|} \overline F(y)dy
&\le& n\overline F(x),
\end{eqnarray}
we obtain statement (i) of the theorem.

Now consider the case (ii) where ${\bf E}\xi\ge 0$.
Take $c>{\bf E}\xi$. Put $\widetilde\xi_i =\xi_i-c$ and
$\widetilde S_n=\widetilde\xi_1+\ldots+\widetilde\xi_n$.
Then ${\bf E}\widetilde\xi={\bf E}\xi-c<0$
and again we can apply Theorem \ref{prop.maxima}.
Thus, there exists a constant $K_1$ such that, for all $x$ and $n$,
\begin{eqnarray*}
\overline{\widetilde F^{*n}}(x) \le
K_1\int_0^{n|{\bf E}\widetilde\xi|}
\overline{\widetilde F}(x+y)dy
\end{eqnarray*}
where $\widetilde F$ in the distribution of $\widetilde\xi$.
Therefore,
\begin{eqnarray*}
{\bf P}\{S_n>x\} = {\bf P}\{\widetilde S_n>x-nc\}
&\le& K_1\int_0^{nc}
\overline{\widetilde F}(x-nc+y)dy\\
&=& K_1\int_0^{nc}
\overline{\widetilde F}(x-y)dy.
\end{eqnarray*}
Since $F\in{\mathcal S}^*$, the distribution $F$ is long-tailed
and, hence, $\overline{\widetilde F}(x)\sim\overline F(x)$
as $x\to\infty$. Then
\begin{eqnarray}\label{upp2}
{\bf P}\{S_n>x\}
&\le& K_2\int_0^{nc}\overline F(x-y)dy,
\end{eqnarray}
for some constant $K_2$ and for all $x\ge 0$.
If $x\ge nc$, then
\begin{eqnarray*}
\int_0^{nc}\overline F(x-y)dy
&\le& \int_0^{nc}\overline F(x-y)
\frac{\overline F(y)}{\overline F(nc)}dy\\
&\le&  \int_0^x \overline F(x-y)
\frac{\overline F(y)}{\overline F(nc)}dy
\le K_3\frac{\overline F(x)}{\overline F(nc)}
\end{eqnarray*}
where
$$
K_3 = \sup_{x\ge 0}
\frac{1}{\overline F(x)}\int_0^x
\overline F(x-y)\overline F(y) dy
$$
is finite, owing to $F\in{\mathcal S}^*$.
If $x<nc$, then
\begin{eqnarray*}
\overline{F^{*n}}(x) &\le& 1
\le \frac{\overline F(x)}{\overline F(nc)}.
\end{eqnarray*}
These two bounds together with (\ref{upp2})
complete the proof of the second assertion
of Theorem \ref{th.asymp}.

From Theorem \ref{th.asymp} and from the dominated
convergence theorem, we deduce the following corollary.

\begin{Corollary}\label{coor}
Tail equivalence {\rm(\ref{eq_lim})} holds
if $F\in{\mathcal S}^*$ and ${\bf E}\xi\ge 0$,
provided that
\begin{eqnarray*}
\sum_{n=1}^\infty \frac{{\bf P}\{\tau =n\}}{\overline F(cn)}
< \infty \quad\mbox{ for some }c>{\bf E}\xi.
\end{eqnarray*}
\end{Corollary}
The latter condition is stronger
than condition (\ref{eq1}), because
\begin{eqnarray*}
\frac{{\bf P}\{\tau>k\}}{\overline F(ck)}
&\le& \sum_{n>k} \frac{{\bf P}\{\tau=n\}}{\overline F(cn)}.
\end{eqnarray*}

Now let us discuss the importance of the condition
$F\in{\mathcal S}^*$ in Theorem \ref{th.asymp}.
The following observation shows the essence
of the difference between two classes of
distributions, ${\mathcal S}$ and ${\mathcal S}^*$,
is the following one. Let a long-tailed
distribution $F$ be absolutely continious with
density $f$. For any function $h(x)>0$,
\begin{eqnarray*}
\int_{h(x)}^{x-h(x)} \overline F(x-y)F(dy)
&=& \int_{h(x)}^{x-h(x)} \overline F(x-y)f(y)dy.
\end{eqnarray*}
Then $F$ is subexponential if and only if
\begin{eqnarray*}
\int_{h(x)}^{x-h(x)} \overline F(x-y)f(y)dy
&=& o(\overline F(x))
\ \mbox{ as }x\to\infty
\end{eqnarray*}
holds for any function $h(x)\to\infty$;
equivalently, if it holds for some function
$h(x)\to\infty$ such that
$\overline F(x-h(x))\sim\overline F(x)$.
On the other hand, $F\in{\mathcal S}^*$
if and only if
\begin{eqnarray*}
\int_{h(x)}^{x-h(x)}
\overline F(x-y)\overline F(y)dy
&=& o(\overline F(x))
\ \mbox{ as }x\to\infty.
\end{eqnarray*}
In typical cases $f(x)=o(\overline F(x))$ and,
hence,
\begin{eqnarray*}
\int_{h(x)}^{x-h(x)} \overline F(x-y)f(y)dy
&=& o\Biggl(\int_{h(x)}^{x-h(x)}
\overline F(x-y)\overline F(y)dy\Biggr)
\ \mbox{ as }x\to\infty.
\end{eqnarray*}
It means that the subexponentiality of $F$
is more likely than $F\in{\mathcal S}^*$.
The latter observation gives the idea how to show
that the condition $F\in{\mathcal S}^*$ in Theorem
\ref{th.asymp}
cannot be extended to the subexponentiality of $F$.

\begin{Theorem}\label{not.n.tails}
There exists a subexponential distribution $F$
on ${\bf R}$ with a negative mean such that
$$
\overline{F^{*n_k}}(x_k)
\ge c\frac{n_k^2}{\ln n_k}\overline F(x_k),
$$
for some $c>0$ and for some sequences $n_k$, $x_k\to\infty$.
\end{Theorem}

The latter theorem yields that, for some distribution
$F\in{\mathcal S}\setminus{\mathcal S}^*$
with negative mean, the first estimate of Theorem
\ref{th.asymp} fails, that is,
$\sup_{n,x}\frac{\overline{F^{*n}}(x)}
{n\overline F(x)}=\infty$.

\proof[Proof \ {\rm of Theorem \ref{not.n.tails}}]
We start with a construction of a specific
subexponential distribution $G$ on the positive half-line.
Put $R_0=0$, $R_1=1$ and $R_{k+1}=e^{R_k}/R_k$ for $k\ge1$.
Since $e^x/x$ is increasing for $x\ge 1$,
the sequence $R_k$ is increasing and
\begin{eqnarray}\label{R.n}
R_k =o(R_{k+1}) \ \mbox{ as }k\to\infty.
\end{eqnarray}
Put $t_k=R_k^2$. Define the hazard function
$R(x)\equiv -\ln\overline G(x)$ as
\begin{eqnarray*}
R(x) = R_k+r_k(x-t_k)
\quad\mbox{ for }x\in(t_k,t_{k+1}],
\end{eqnarray*}
where
\begin{eqnarray}\label{rn.def}
r_k &=& \frac{R_{k+1}-R_k}{t_{k+1}-t_k}
=\frac{1}{R_{k+1}+R_k}
\sim \frac{1}{R_{k+1}}
\end{eqnarray}
by (\ref{R.n}).
In other words, the hazard rate $r(x)=R'(x)$
is defined as $r(x)=r_k$ for $x\in(t_k,t_{k+1}]$,
where $r_k$ is given by (\ref{rn.def}).
By the construction, we have
$\overline G(t_k) = e^{-\sqrt{t_k}}$,
so that at points $t_k$ the tail of $G$ behaves
like the Weibull tail with parameter $1/2$.
Between these points the tail decays
exponentially with indexes $r_k$.

We prove now that $G$ has finite mean
and is subexponential. Since by (\ref{rn.def})
\begin{eqnarray*}
\int_{t_k}^{t_{k+1}} e^{-R(y)} dy
&=& r_k^{-1} (e^{-R_k}-e^{-R_{k+1}})\nonumber\\
&\sim& r_k^{-1} e^{-R_k} \sim R_{k+1} e^{-R_k}=1/R_k,
\end{eqnarray*}
the mean of $G$,
\begin{eqnarray*}
\int_0^\infty \overline G(y)dy
=\sum_{k=0}^\infty\int_{t_k}^{t_{k+1}}\overline G(y)dy,
\end{eqnarray*}
is finite.

It follows from the definition that $r(x)$
decreases to 0. Then we can apply
Pitman's criterion [\ref{Pitman}] which says that
$G$ is subexponential if the function
$e^{yr(y)-R(y)}r(y)$ is integrable over $[0,\infty)$.
In order to estimate the integral of this function, put
\begin{eqnarray*}
I_k = \int_{t_k}^{t_{k+1}} e^{yr(y)-R(y)}r(y)dy.
\end{eqnarray*}
Then
\begin{eqnarray*}
I_k = r_k\int_{t_k}^{t_{k+1}} e^{yr_k-(R_k+r_k(y-t_k))}dy
\le r_k e^{-R_k+r_kt_k}t_{k+1}.
\end{eqnarray*}
Since
\begin{eqnarray}\label{rn.tn+1}
r_kt_{k+1}=r_kR^2_{k+1} \sim R_{k+1}
\end{eqnarray}
by (\ref{rn.def}) and
\begin{eqnarray*}
r_kt_k=r_kR^2_k \sim R_k^2/R_{k+1}
=R_k^3 e^{-R_k} \to 0,
\end{eqnarray*}
we get $I_k \le 2R_{k+1} e^{-R_k}\sim 2/R_k$
for $k$ sufficiently large.
Therefore,
\begin{eqnarray*}
\int_0^\infty e^{yr(y)-R(y)}r(y)dy
= \sum_{k=0}^\infty I_k <\infty,
\end{eqnarray*}
and $G$ is indeed subexponential.

In the sequel we need to know the asymptotic
behaviour of the following internal part
of the convolution integral at point $t_k$:
\begin{eqnarray*}
J_k=\int_{t_k/4}^{3t_k/4}
\overline G(t_k-y)G(dy)=
\int_{t_k/4}^{3t_k/4}
e^{-R(t_k-y)}e^{-R(y)}r(y)dy.
\end{eqnarray*}
Owing to (\ref{R.n}), $t_{k-1}=o(t_k)$.
Thus, $(t_k/4,3t_k/4]\subset(t_{k-1},t_k-t_{k-1}]$
for all sufficiently large $k$. For those values
of $k$, we have
\begin{eqnarray*}
J_k &=& \overline G(t_k)\int_{t_k/4}^{3t_k/4}
e^{-(-r_{k-1}y)}e^{-(R_{k-1}+r_{k-1}(y-t_{k-1}))}r_{k-1}dy\\
&\ge& \overline G(t_k)(t_k/2)e^{-R_{k-1}}r_{k-1}.
\end{eqnarray*}
Applying (\ref{rn.tn+1}) and the equality
$e^{R_{k-1}}=R_kR_{k-1}$,
we obtain, for all sufficiently large $k$,
\begin{eqnarray}\label{inter.int.4}
J_k \ge \overline G(t_k)e^{-R_{k-1}}R_k/3
= \overline G(t_k)/3R_{k-1}.
\end{eqnarray}

Let $\eta_1$, $\eta_2$, \ldots\ be independent
random variables with common distribution $G$
and put $T_n=\eta_1+\ldots+\eta_n$.
For any $n$, we have
\begin{eqnarray*}
{\bf P}\{T_n>x\} &\ge& \sum_{1\le i<j\le n}
{\bf P}\{T_n>x,\eta_i>n,\eta_j>n,
\eta_l\le n\mbox{ for all }l\not=i,j\}\\
&=& \frac{n(n-1)}{2}
{\bf P}\{T_n>x,\eta_1>n,\eta_2>n,
\eta_3\le n,\ldots,\eta_n\le n\}.
\end{eqnarray*}
Since $\eta$'s are positive,
the latter probability is not smaller than
\begin{eqnarray*}
{\bf P}\{\eta_1+\eta_2>x,\eta_1>n,\eta_2>n\}
{\bf P}\{\eta_3\le n,\ldots,\eta_n\le n\}.
\end{eqnarray*}
The mean of $\eta$ is finite, thus
$\overline G(n)=o(1/n)$ as $n\to\infty$ and
\begin{eqnarray*}
{\bf P}\{\eta_3\le n,\ldots,\eta_n\le n\}
= (1-\overline G(n))^{n-2} \to 1.
\end{eqnarray*}
Putting altogether, we get, for all sufficiently
large $n$, the following estimate from below
\begin{eqnarray}\label{T.k.lower}
{\bf P}\{T_n>x\}
&\ge& \frac{n^2}{3}{\bf P}
\{\eta_1+\eta_2>x,\eta_1>n,\eta_2>n\}.
\end{eqnarray}

Now take $n=n_k=[\sqrt{t_k}]=[R_k]$.
Then, for all sufficiently large $k$
(at least for those $k$ where $n_k<t_k/4$),
\begin{eqnarray*}
{\bf P}\{\eta_1+\eta_2>t_k,\eta_1>n_k,\eta_2>n_k\}
\ge J_k.
\end{eqnarray*}
Therefore, by (\ref{T.k.lower}) and (\ref{inter.int.4}),
for all sufficiently large $k$,
\begin{eqnarray*}
{\bf P}\{T_{n_k}>t_k\}
\ge n_k^2 \overline G(t_k)/9R_{k-1}
\sim n_k^2\overline G(t_k)/9\ln n_k,
\end{eqnarray*}
due to $R_{k-1}\sim\ln R_k\sim\ln n_k$.

Denote $b={\bf E}\eta_1$. Put $\xi_i=\eta_i-2b$,
then $\xi$'s have negative mean and $S_n=T_n-2nb$.
Denote by $F$ the distribution of $\xi_1$;
it is subexponential because $G$ is.

Take $x=x_k=t_k-2n_kb$, so that $x_k\sim n_k^2$.
By the latter inequality we have
\begin{eqnarray*}
{\bf P}\{S_{n_k}>x_k\} = {\bf P}\{T_{n_k}>t_k\}
\ge n_k^2\overline G(t_k)/10\ln n_k.
\end{eqnarray*}
Note also that
\begin{eqnarray*}
\overline F(x_k)=\overline G(t_k-2n_kb)
=\overline G(t_k)e^{r_{k-1}2n_kb}\le
\overline G(t_k) e^{2b}
\end{eqnarray*}
because $r_{k-1}n_k\le r_{k-1}R_k\le 1$ by (\ref{rn.def}).
Therefore, the inequality
\begin{eqnarray*}
{\bf P}\{S_{n_k}>x_k\}
\ge n_k^2\overline F(x_k) e^{-2b}/10\ln n_k
\end{eqnarray*}
holds which yields the conclusion of the theorem.

The subexponential distribution $G$
constructed in the latter proof cannot belong
to the class ${\mathcal S}^*$ because otherwise
the theorem conclusion fails, as follows
from Theorem \ref{th.asymp}. The fact that
$G\not\in{\mathcal S}^*$ can also be proved directly.
Kl\"uppelberg's criterion [\ref{Kl}] states
that $G\in{\mathcal S}^*$ if and only if
\begin{eqnarray*}
\int_0^x e^{yr(x)-R(y)}dy
&\to& \int_0^\infty\overline G(y)dy
\quad \mbox{ as }x\to\infty.
\end{eqnarray*}
In our construction,
\begin{eqnarray*}
\int_0^{t_k-0} e^{yr(t_k-0)-R(y)}dy
&\ge& \int_{t_{k-1}}^{t_k} e^{yr_{k-1}-R(y)}dy\\
&\ge& (t_k-t_{k-1})e^{-R_{k-1}}\\
&\sim& R_k^2 e^{-R_{k-1}}
=e^{R_{k-1}}/R_{k-1}^2\to\infty
\end{eqnarray*}
as $k\to\infty$. Hence, $G\not\in{\mathcal S}^*$.

\mysection{On the asymptotics
${\bf P}\{S_n>x\}\sim n\overline F(x)$}\label{nFx}

As before, we assume ${\bf E}\xi$ to be finite.
Then, by the Strong Law of Large Numbers,
\begin{eqnarray}\label{stoch.bound}
{\bf P}\{S_n>-An\} \to 1
\ \mbox{ as }A\to\infty\mbox{ uniformly in }n\ge1,
\end{eqnarray}
and by the Chebyshev's inequality
\begin{eqnarray}\label{bound.for.one}
{\bf P}\{\xi_1>An\} \le {\bf E}|\xi_1|/An
\ \mbox{ for all }A>0\mbox{ and }n\ge1.
\end{eqnarray}

\begin{Theorem}\label{equiv.nF}
Let $F\in{\mathcal S}^*$ and let an increasing
function $h(x)>0$ be such that
$\overline F(x\pm h(x))\sim\overline F(x)$.
Then ${\bf P}\{S_n>x\}\sim n\overline F(x)$
as $x\to\infty$ uniformly in $n\le h(x)$.
\end{Theorem}

\proof[Proof \ {\rm of the lower bound is similar to
that in {[\ref{DDS}, Section 4]}}]
Fix $A>0$. We use the following inequalities:
\begin{eqnarray*}
{\bf P}\{S_n>x\} &\ge& \sum_{i=1}^n
{\bf P}\{S_n>x,\xi_i>x+An,\xi_j\le An
\mbox{ for all }j\not= i\}\\
&\ge& n{\bf P}\{S_n-\xi_1>-An,
\xi_1>x+An,\xi_2\le An,\ldots,\xi_n\le An\}\\
&=& n\overline F(x+An) {\bf P}\{S_{n-1}>-An,
\xi_1\le An,\ldots,\xi_{n-1}\le An\}.
\end{eqnarray*}
We have $\overline F(x+An)\sim\overline F(x)$ as
$x\to\infty$ uniformly in $n\le h(x)$.
Taking also into account that
\begin{eqnarray*}
{\bf P}\{S_{n-1}>-An,
\xi_1\le An,\ldots,\xi_{n-1}\le An\}
&\ge& {\bf P}\{S_{n-1}>-An\}-(n-1){\bf P}\{\xi_1>An\},
\end{eqnarray*}
we get, for any fixed $A>0$,
\begin{eqnarray*}
\liminf_{x\to\infty}\inf_{n\le h(x)}
\frac{{\bf P}\{S_n>x\}}{n\overline F(x)}
&\ge& \inf_n
({\bf P}\{S_{n-1}>-An\}-(n-1){\bf P}\{\xi_1>An\}).
\end{eqnarray*}
Since the infimum on the right goes to $1$ as $A\to\infty$
owing to (\ref{stoch.bound}) and (\ref{bound.for.one}),
we arrive at the following lower bound
\begin{eqnarray*}
\liminf_{x\to\infty}\inf_{n\le h(x)}
\frac{{\bf P}\{S_n>x\}}{n\overline F(x)}
&\ge& 1.
\end{eqnarray*}

To prove the upper bound, we apply
Theorem \ref{prop.maxima} to random variables
$\widetilde\xi_i=\xi_i-{\bf E}\xi_1-1$ with negative
mean ${\bf E}\widetilde\xi_i=-1$ and to
$\widetilde S_n=S_n-n({\bf E}\xi_1+1)$. Thus,
\begin{eqnarray*}
{\bf P}\{S_n>x\} &=& {\bf P}\{\widetilde S_n>x-n({\bf E}\xi_1+1)\}\\
&\le& (1+o(1))\int_{x-n({\bf E}\xi_1+1)}^{x-n{\bf E}\xi_1}
\overline{\widetilde F}(x+u)du\\
&\le& (1+o(1))n\overline{\widetilde F}(x-n({\bf E}\xi_1+1))
\end{eqnarray*}
as $x\to\infty$ where $\widetilde F$ is
the  distribution of $\widetilde\xi$.
If $n\le h(x)$ then
$\overline{\widetilde F}(x-n({\bf E}\xi_1+1))
\sim\overline F(x)$
as $x\to\infty$ and the proof is complete.

The range $n\le h(x)$ is usually more narrow
than one could expect. Say, for the regularly
varying distributions (more generally,
for the intermediate regularly varying,
see the definition in Section \ref{tailcom}) we can
take $h(x)=o(x)$. Then we get the range
$n=o(x)$ while the standard (if the mean is zero
and the second moment is finite) range is $x^2>cn\ln n$;
in the class of distributions with finite mean,
the relation ${\bf P}\{S_n>x\}\sim n\overline F(x)$ holds
in the range $x>({\bf E}\xi+\varepsilon)n$,
$\varepsilon>0$, see S. Nagaev [\ref{SN1981}].
The advantage of the result in Theorem \ref{equiv.nF}
is its simplicity and universality since it is
valid for all distributions from ${\mathcal S}^*$
without any further moment or regularity assumptions,
compare with a series of results in
[\ref{BB}, \ref{BM}, \ref{DDS}] where the hazard
rate is assumed to be sufficiently smooth.

As follows from [\ref{DDS}], if the mean is zero
and the second moment is finite, then the right range
should be $n\le h^2(x)$, roughly speaking.
Our technique allows to prove the lower bound
for this range.

\begin{Theorem}\label{tm.lower.nF}
Let ${\bf E}\xi=0$ and ${\bf E}\xi^2<\infty$.
Let $F$ be a long-tailed distribution
and let an increasing function $h(x)>0$ be such that
$\overline F(x\pm h(x))\sim\overline F(x)$.
Then ${\bf P}\{S_n>x\}\ge(1+o(1))n\overline F(x)$
as $x\to\infty$ uniformly in $n\le h^2(x)$.
\end{Theorem}

\proof
Fix $A>0$. By the Chebyshev's inequality,
\begin{eqnarray}\label{stoch.bound2}
{\bf P}\{\xi_1>A\sqrt n\}\le {\bf E}\xi^2/A^2n
\ \mbox{ and }\
{\bf P}\{S_n>-A\sqrt n\} \ge 1-{\bf E}\xi^2/A^2.
\end{eqnarray}
In this proof we use a slightly different inequality
than in the previous theorem:
\begin{eqnarray*}
{\bf P}\{S_n>x\} &\ge& \sum_{i=1}^n
{\bf P}\{S_n>x,\xi_i>x+A\sqrt n,\xi_j\le A\sqrt n
\mbox{ for all }j\not= i\}\\
&\ge& n{\bf P}\{S_n-\xi_1>-A\sqrt n,
\xi_1>x+A\sqrt n,\xi_2\le A\sqrt n,
\ldots,\xi_n\le A\sqrt n\}\\
&=& n\overline F(x+A\sqrt n)
{\bf P}\{S_{n-1}>-A\sqrt n,\xi_1\le A\sqrt n,
\ldots,\xi_{n-1}\le A\sqrt n\}.
\end{eqnarray*}
Since $n\le h^2(x)$,
$\overline F(x+A\sqrt n)\sim\overline F(x)$
as $x\to\infty$. Applying (\ref{stoch.bound2}), we get
\begin{eqnarray*}
\lefteqn{{\bf P}\{S_{n-1}>-A\sqrt n,\xi_1\le A\sqrt n,
\ldots,\xi_{n-1}\le A\sqrt n\}}\\
&\ge& {\bf P}\{S_{n-1}>-A\sqrt n\}-
(n-1){\bf P}\{\xi_1>A\sqrt n\}\\
&\ge& 1-2{\bf E}\xi^2/A^2\to 1
\ \mbox{ as }A\to\infty.
\end{eqnarray*}
Now the lower bound for ${\bf P}\{S_n>x\}$ follows.

\mysection{Proof of Theorem \ref{thm1}}\label{pth1}

Since $\tau$ is independent of $\xi$'s,
we can use the following decomposition:
$$
{\bf P}\{S_\tau>x\} = \sum_{n=0}^\infty
{\bf P}\{\tau=n\}\overline{F^{*n}}(x).
$$
By the subexponentiality, here the $n$th term is
equivalent to $n{\bf P}\{\tau=n\}\overline F(x)$
as $x\to\infty$. In particular, by Fatou's lemma,
\begin{eqnarray}\label{lim.inf.th}
\liminf_{x\to\infty} \frac{{\bf P}\{S_\tau>x\}}{\overline F(x)}
\ge \sum_{n=0}^\infty n{\bf P}\{\tau=n\}={\bf E}\tau,
\end{eqnarray}
without any condition on the sign of ${\bf E}\xi$.
In the case of negative mean, the $n$th term
is bounded from above by $n\overline F(x)$,
see (\ref{ine.n}).
Then the dominated convergence for series
yields statement (i) of the theorem.

Now turn to the proof of statement (ii)
where ${\bf E}\xi\ge 0$.
Since $S_\tau \le M_\tau$,
it follows from (\ref{lim.inf.th})
that it is sufficient to prove that
\begin{equation}\label{two}
{\bf P}\{M_\tau>x\} \sim
{\bf E}\tau\overline F(x)
\ \mbox{ as }x\to\infty.
\end{equation}
To prove the latter relation,
we start with the following representation:
for any $N$,
\begin{eqnarray}\label{prob=sum}
{\bf P}\{M_\tau>x\}
&=& {\bf P}\{M_\tau>x, \tau\le N\}
+{\bf P}\{M_\tau>x, \tau\in(N,x/c]\}
+{\bf P}\{M_\tau>x, c\tau>x\}\nonumber\\
&\equiv& P_1+P_2+P_3.
\end{eqnarray}

Since any ${\mathcal S}^*$-distribution is subexponential
and $S_n\le M_n\le \xi_1^++\ldots+\xi_n^+$,
$$
{\bf P}\{M_n>x\}\sim n\overline F(x)
$$
as $x\to\infty$, for any $n$. Thus, for any fixed $N$,
\begin{eqnarray*}
{\bf P}\{M_\tau>x, \tau\le N\}
= \sum_{n=1}^N {\bf P}\{\tau=n\}{\bf P}\{M_n>x\}
\sim {\bf E}\{\tau;\tau\le N\}\overline F(x)
\end{eqnarray*}
as $x\to\infty$ which implies the existence of
an increasing function $N(x)\to\infty$ such that
\begin{eqnarray}\label{first.term}
P_1={\bf P}\{M_\tau>x, \tau\le N(x)\}
\sim {\bf E}\tau\overline F(x).
\end{eqnarray}

In what follows, we use representation
(\ref{prob=sum}) with $N(x)$ in place of $N$.
We further estimate the second term on
the right side in (\ref{prob=sum}).
Let $\varepsilon=(c-{\bf E}\xi)/2>0$
and let $b=({\bf E}\xi+c)/2$.
Consider $\widetilde \xi_n =\xi_n-b$,
$\widetilde S_n=\widetilde \xi_1+\ldots+\widetilde\xi_n$ and
$\widetilde M_n=\max (\widetilde S_1,\ldots,\widetilde S_n)$.
Then ${\bf E}\widetilde\xi=-\varepsilon<0$
and we can apply Theorem~\ref{prop.maxima}.
Taking into account that
$M_n\le\widetilde M_n+bn$,
we obtain that there exists $K$ such that,
for all $x$ and $n$,
\begin{eqnarray*}
{\bf P}\{M_n>x\}
&\le& {\bf P}\{\widetilde M_n>x-bn\}\\
&\le& K\int_0^{n\varepsilon}\overline{\widetilde F}(x-nb+y)dy\\
&\le& K\int_0^{n\varepsilon}\overline F(x-nb+y)dy.
\end{eqnarray*}
Hence,
\begin{eqnarray*}
P_2={\bf P}\{M_\tau>x, \tau\in(N(x),x/c]\}
&\le& K\sum_{n=N(x)}^{[x/c]}
{\bf P}\{\tau=n\} \int_0^{n\varepsilon}
\overline F(x-nb+y)dy.
\end{eqnarray*}
Since $b-\varepsilon={\bf E}\xi$,
\begin{eqnarray*}
\int_0^{n\varepsilon}\overline F(x-nb+y)dy
&=& \int_{n{\bf E}\xi}^{nb}\overline F(x-y)dy.
\end{eqnarray*}
Then
\begin{eqnarray}\label{P_2}
P_2 &\le& K\int_{N(x){\bf E}\xi}^{b[x/c]} \overline F(x-y)dy
\sum_{n=\max(N(x),[y/b]+1)}^{[x/c]} {\bf P}\{\tau=n\}
\nonumber\\
&\le& K\int_{N(x){\bf E}\xi}^{bx/c} \overline F(x-y)
{\bf P}\{\tau> y/b\}dy\nonumber\\
&\le& K\int_{N(x){\bf E}\xi}^{bx/c} \overline F(x-y)
{\bf P}\{\tau> y/c\}dy,
\end{eqnarray}
because $b<c$.
By condition (\ref{eq1}),
${\bf P}\{\tau>y/c\} \le K_1\overline F(y)$,
for some $K_{1}$ and all $y$.
Therefore, the inequality
\begin{eqnarray}\label{second.term}
P_2 \le KK_1\int_{N(x){\bf E}\xi}^{bx/c}
\overline F(x-y)\overline F(y)dy
= o(\overline F(x))
\ \mbox{ as }x\to\infty
\end{eqnarray}
follows from $b/c<1$ and from $F\in{\mathcal S}^*$.
Indeed, for any ${\mathcal S}^*$-distribution,
\begin{equation}\label{hstar}
\int_{h(x)}^{x-h(x)}  \overline F(x-y)\overline F(y)dy
= o(\overline F(x))\ \mbox{ as }x\to\infty,
\end{equation}
for any function $h(x)\to\infty$ such that $h(x)\le x/2$
(see, e.g., [\ref{Kl}]).

Now we estimate the third term on the right in
(\ref{prob=sum}) using condition (\ref{eq1}):
\begin{eqnarray}\label{third.term}
P_3 \le {\bf P}\{c\tau>x\} = o(\overline F(x))
\ \mbox{ as }x\to\infty.
\end{eqnarray}
Altogether relations (\ref{first.term}),
(\ref{second.term}), and (\ref{third.term})
complete the proof of Theorem \ref{thm1}.

Now we provide an example where
$$
\frac{{\bf P}\{S_\tau>x\}}{\overline F(x)} \to \infty
$$
given that condition (\ref{eq1}) is satisfied only
with $c={\bf E}\xi>0$ and not with any bigger $c$.
Assume that $F$ is a Weibull distribution on the positive half line
with parameter $\beta\in(1/2,1)$, that is
$\overline F(x) = e^{-x^\beta}$. Let $\tau$
have a distribution such that
${\bf P}\{c\tau>x\} \sim x^{-1} e^{-x^\beta}$
as $x\to\infty$.
Write down the following lower bound:
\begin{eqnarray*}
{\bf P}\{S_\tau>x\}
&\ge& {\bf P}\{S_\tau>x|c\tau>x-\sqrt x\}
{\bf P}\{c\tau>x-\sqrt x\}.
\end{eqnarray*}
By the Central Limit Theorem,
\begin{eqnarray*}
\delta \equiv \liminf_{x\to\infty}
{\bf P}\{S_\tau>x|c\tau>x-\sqrt x\}
\ge \liminf_{x\to\infty}
{\bf P}\{S_{[(x-\sqrt x)/c]}>x\}>0.
\end{eqnarray*}
Hence,
\begin{eqnarray*}
\liminf_{x\to\infty}
\frac{{\bf P}\{S_\tau>x\}}{\overline F(x)}
&\ge& \delta \liminf_{x\to\infty}
\frac{{\bf P}\{c\tau>x-\sqrt x\}}{\overline F(x)}\\
&=& \delta \liminf_{x\to\infty}
\frac{e^{x^\beta-(x-\sqrt x)^\beta}}
{x-\sqrt x}=\infty,
\end{eqnarray*}
because $\beta>1/2$.

We conclude this section by an example
showing that the conclusion of Theorem \ref{thm1}
cannot hold for all subexponential distributions.
Indeed, take $F$ with negative mean as described
in Theorem \ref{not.n.tails}. Without loss
of generality we assume that the series
$\sum_k n_k^{-1}\ln n_k$ converge.
Consider $\tau$ taking values
$n_k$ with probabilities $c\ln^2 n_k/n_k^2$,
here $c$ is the normalising constant.
Then $\tau$ has a finite mean, but
\begin{eqnarray*}
{\bf P}\{S_\tau>x_k\}\ge
{\bf P}\{S_{n_k}>x_k\}{\bf P}\{\tau=n_k\}
\ge c\frac{n_k^2}{\ln n_k}\overline F(x_k)
\frac{\ln^2 n_k}{n_k^2},
\end{eqnarray*}
so that, as $k\to\infty$,
\begin{eqnarray*}
\frac{{\bf P}\{S_\tau>x_k\}}
{\overline F(x_k)} \to \infty.
\end{eqnarray*}

\mysection{The case where $\xi$ and $\tau$
may be tail-comparable}\label{tailcom}

In this section we do not assume condition
(\ref{eq1}) to hold, such a situation is
of particular importance for branching processes.
To start with, we define two important classes
of distributions.

A distribution $F$ is called {\it dominated varying}
if there exists $c$ such that
$\overline F(x) \le c\overline F(2x)$ for all $x$.
It is known that any long-tailed and dominated varying
distribution with a finite mean
belongs to the class ${\mathcal S}^*$,
see [\ref{Kl}].

We say that a distribution $G$ is {\it intermediate
regularly varying} (at infinity) if
\begin{eqnarray}\label{i.r.v}
\lim_{\varepsilon\downarrow 0}\limsup_{x\to\infty}
\frac{\overline G((1-\varepsilon)x)}{\overline G(x)}
= 1.
\end{eqnarray}
In particular, any regularly varying at infinity
distribution satisfies the latter relation.
Any intermediate regularly varying distribution
is long-tailed and dominated varying;
in particular, it belongs to the class ${\mathcal S}^*$,
provided its mean is finite.

\begin{Theorem}\label{thm2}
Let $F\in{\mathcal S}^*$, ${\bf E}\xi>0$, and
\begin{eqnarray}\label{F.le.tau}
\overline F(x)=O({\bf P}\{\tau>x\})
\ \mbox{ as }x\to\infty.
\end{eqnarray}
If the distribution of $\tau$ is
intermediate regularly varying, then
\begin{eqnarray}\label{comparable}
{\bf P}\{S_\tau>x\} \sim {\bf P}\{M_\tau>x\} \sim
{\bf E}\tau\overline F(x)+{\bf P}\{\tau>x/{\bf E}\xi\}
\ \mbox{ as }x\to\infty.
\end{eqnarray}
\end{Theorem}

We strongly believe that the statement of the theorem
stays valid in a more general setting where
the distribution of $\tau$ is assumed to be
{\it square-root insensitive}, that is
${\bf P}\{\tau>x\pm \sqrt x\}\sim{\bf P}\{\tau>x\}$,
and the variance of $\xi$ is finite. Probably,
some further minor regularity assumptions are required.
For example, the Weibull distribution
$\overline F(x)=e^{-x^\beta}$ with parameter
$\beta<1/2$ is square-root insensitive.
For distribution which is {\it not}
square-root insensitive, the asymptotics are
different and more complicated.

\proof[Proof of Theorem \ref{thm2}]
By (\ref{i.r.v}), for any fixed $\delta>0$,
we can choose $a<{\bf E}\xi$ and $c>{\bf E}\xi$
sufficiently close to ${\bf E}\xi$ such that
\begin{eqnarray*}
1-\delta/2 \le
\liminf_{x\to\infty}\frac{{\bf P}\{a\tau>x\}}
{{\bf P}\{\tau>x/{\bf E}\xi\}}
\le \limsup_{x\to\infty}\frac{{\bf P}\{c\tau>x\}}
{{\bf P}\{\tau>x/{\bf E}\xi\}}
\le 1+\delta/2.
\end{eqnarray*}
Then, due to $S_\tau \le M_\tau$, it is sufficient
to prove the following lower bound for the sum
\begin{eqnarray}\label{l.b.S}
{\bf P}\{S_\tau>x\} &\ge&
({\bf E}\tau+o(1))\overline F(x)
+(1+o(1)){\bf P}\{\tau>x/a\}.
\end{eqnarray}
and the upper bound for the maximum
\begin{eqnarray}\label{u.b.M}
{\bf P}\{M_\tau>x\} &\le&
({\bf E}\tau+o(1))\overline F(x)
+(1+o(1)){\bf P}\{\tau>x/c\}
\ \mbox{ as }x\to\infty.
\end{eqnarray}

We have
\begin{eqnarray*}
{\bf P}\{S_\tau>x\}
&=& {\bf P}\{S_\tau>x, \tau\le x/a\}
+{\bf P}\{S_\tau>x, \tau>x/a\}.
\end{eqnarray*}
Since $a<{\bf E}\xi$,
${\bf P}\{S_\tau>x|\tau>x/a\}\to1$
as $x\to\infty$, by the Law of Large Numbers.
Now the standard arguments lead to (\ref{l.b.S}).

To prove the upper bound, we use a representation
similar to (\ref{prob=sum}) (see the previous proof):
\begin{eqnarray*}
{\bf P}\{M_\tau>x\}
&=& {\bf P}\{M_\tau>x, \tau\le N(x)\}
+{\bf P}\{M_\tau>x, \tau\in(N(x),x/c]\}
+{\bf P}\{M_\tau>x, c\tau>x\}\\
&\equiv& P_1+P_2+P_3.
\end{eqnarray*}
The first summand $P_1$ can be treated as ealier.
The second summand $P_2$ can be estimated as follows:
if condition (\ref{F.le.tau}) holds then,
by estimate (\ref{P_2}),
\begin{eqnarray*}
P_2 &\le& KK_2\int_{N(x){\bf E}\xi}^{bx/c}
{\bf P}\{\tau> x-y\}
{\bf P}\{\tau> y\}dy,
\end{eqnarray*}
for some $K_2$. Since the distribution of
$\tau$ is intermediate regularly varying and,
therefore, belongs to ${\mathcal S}^*$,
\begin{eqnarray*}
P_2 = o({\bf P}\{\tau>x\}).
\end{eqnarray*}
Taking into account also that $P_3\le{\bf P}\{c\tau>x\}$,
we finally get
\begin{eqnarray*}
{\bf P}\{M_\tau>x\} &\le&
({\bf E}\tau+o(1))\overline F(x)
+{\bf P}\{\tau>x/c\}+o({\bf P}\{\tau>x\})
\ \mbox{ as }x\to\infty.
\end{eqnarray*}
Since the distribution of $\tau$ is (in particular)
dominated varying, ${\bf P}\{\tau>x\}=O({\bf P}\{\tau>x/c\})$.
Therefore, (\ref{u.b.M}) is proved and
the conclusion of Theorem \ref{thm2} follows.

\begin{Theorem}\label{thm2.D}
Let ${\bf E}\xi>0$ and let
$\tau$ have an intermediate regularly
varying distribution. If the distribution
$F$ is long-tailed and dominated varying,
then {\rm(\ref{comparable})} holds.
\end{Theorem}

A particular corollary is that if both $\xi$
and $\tau$ have regularly varying tail distributions,
then asymptotics (\ref{comparable}) hold;
this result was proved by Stam [\ref{Stam}, Theorems 1.3 and 1.4]
for positive $\xi$ and by A. Borovkov and K. Borovkov
[\ref{BB}, Section 7.1] for signed $\xi$.
Also, Theorems \ref{thm2} and \ref{thm2.D}
generalise and improve Theorem 1.3 of
Ale\v skevi\v cen\.e {\it et al.} [\ref{ALS}].

\proof[Proof of Theorem \ref{thm2.D}]
It follows the lines of the previous proof,
and only the term $P_2$ needs a different estimation.
From bound (\ref{P_2}), we get
\begin{eqnarray*}
P_2 &\le& K\overline F(x-bx/c)
\int_{N(x){\bf E}\xi}^{bx/c} {\bf P}\{c\tau> y\}dy.
\end{eqnarray*}
Since $F$ is dominated varying,
$\overline F(x-bx/c)=O(\overline F(x))$ as $x\to\infty$.
Therefore, $P_2=o(\overline F(x))$
and the proof is complete.

\mysection{Applications to the branching processes}\label{branch}

A Galton--Watson process is a stochastic process
$\{X_n\}$ which evolves according to the recurrence
formula $X_0=1$ and
\begin{eqnarray*}
X_{n+1} &=&
\sum_{j=1}^{X_n}\xi_j^{(n+1)},
\end{eqnarray*}
where $\{\xi_j^{(n)}\}$ is a family of independent
identically distributed non-negative integer-valued
random variables with a finite mean,
and their common distribution does not depend on $n$.
Here $X_n$ is the number of items in the $n$th generation.
Taking into account that
any intermediate regularly varying
distribution with finite mean
belongs to the class ${\mathcal S}^*$,
we obtain the following application of Theorem \ref{thm2}
to the branching process:

\begin{Corollary}\label{co:1}
Let the common distribution of $\xi$'s be
intermediate regularly varying.
Then, as $x\to\infty$,
\begin{eqnarray*}
{\bf P}\{X_2>x\} \sim
{\bf E}\xi{\bf P}\{\xi>x\}+{\bf P}\{\xi>x/{\bf E}\xi\}.
\end{eqnarray*}
\end{Corollary}

In particular, if the branching process is
{\it critical}, i.e. if ${\bf E}\xi=1$, then
\begin{eqnarray*}
{\bf P}\{X_2>x\} \sim 2{\bf P}\{\xi>x\}
\ \mbox{ as }x\to\infty.
\end{eqnarray*}
More generally, by induction arguments,
the tail of the distribution of the number
of items in the $n$th generation is
asymptotically equivalent to $n{\bf P}\{\xi>x\}$.
A similar result (for critical process) was obtained in
[\ref{Wachtel}, Theorem 2] in the case of regularly varying
distribution of $\xi$'s and for possibly growing $n$.

\mysection{Equivalences in the case where
a counting random variable $\tau$
may depend on $\xi$'s}\label{taudep}

We continue to assume that random variables $\{\xi_n\}$
are independent and identically distributed.
For any family $\Xi$ of random variables,
denote by $\sigma(\Xi)$ the $\sigma$-algebra
generated by $\Xi$.
Traditionally, a counting random variable $\tau$ is called
a {\it stopping time} for a sequence $\{\xi_n\}$ if
$\{\tau\le n\}\in\sigma(\xi_1,\ldots,\xi_n)$ for all $n$.

We say that a counting random variable $\tau$
{\it does not depend on the future of the sequence}
$\{\xi_n\}$ if the family
$(\xi_1,\ldots,\xi_n,{\bf I}\{\tau\le n\})$
does not depend on $(\xi_j,j\ge n+1)$ for all $n$.
Dependence of this type goes back to Kolmogorov
and Prokhorov [\ref{KP}] who proved Wald's
identity under the condition that the event
$\{\tau\le n\}$ does not depend on
$\xi_j$ for all $n\ge 1$ and $j\ge n+1$.

Provided independence of $\xi$'s, any stopping
time $\tau$ does not depend on the future
of the sequence $\{\xi_n\}$.
If a counting random variable $\tau$ does not
depend on $\xi$'s, then it does not depend
on the future of the sequence $\{\xi_n\}$.

Let ${\mathcal F}_n$ be a filtration of $\sigma$-algebras.
A counting random variable $\tau$ is called
a stopping time for this filtration
if $\{\tau\le n\}\in{\mathcal F}_n$ for all $n$.
In this terminology, $\tau$ is a stopping time
for a sequence $\{\xi_n\}$ if and only if $\tau$
is a stopping time for the natural filtration
${\mathcal F}_n=\sigma(\xi_1,\ldots,\xi_n)$.

Consider a special filtration
${\mathcal F}_n=\sigma(\xi_k,{\bf I}\{\tau=k\},k\le n)$.
Then $\tau$ is a stopping time for this filtration.
In addition, $\tau$ does not depend on the future
of the sequence $\{\xi_n\}$ if and only if
$(\xi_j,j\ge n+1)$ does not depend
on ${\mathcal F}_n$ for all $n$.

We start with a result for a bounded counting
stopping time (recall that a random variable is {\it bounded}
if its distribution has a bounded support).

\begin{Theorem}\label{thm.stopping}
Let $\xi$ have a subexponential distribution $F$
on ${\bf R}$ (we do not assume finite mean),
and let the counting variable $\tau$ do not
depend on the future. If $\tau$ is bounded, then
${\bf P}\{S_\tau>x\}\sim{\bf E}\tau\overline F(x)$
as $x\to\infty$.
\end{Theorem}

Similar result for $M_\tau$ may be found
in [\ref{FPZ}, Theorem 1].
Note that one cannot expect the latter asymptotics
to hold for any $\tau$ with unbounded support,
which may depend on $\{\xi_n\}$ -- even for
a stopping time. Indeed, consider a stopping
time $\tau=\min\{n:S_n\le0\}$.
If ${\bf E}\xi<0$ then ${\bf E}\tau$ is finite but
${\bf P}\{S_\tau>x\} =0$ for any $x>0$.

\proof
We adopt the corresponding proof from
Greenwood [\ref{Greenwood}] where a stopping
time and regularly varying tails were considered.
Let $N$ be such that ${\bf P}\{\tau\le N\}=1$.
The starting point of the proof
is the following representation:
\begin{eqnarray*}
{\bf P}\{S_\tau>x\} &=&
\sum_{n=1}^N ({\bf P}\{S_n>x,\tau\ge n\}
-{\bf P}\{S_n>x,\tau\ge n+1\})\\
&=& {\bf P}\{S_1>x,\tau\ge 1\}
+\sum_{n=2}^N ({\bf P}\{S_n>x,\tau\ge n\}
-{\bf P}\{S_{n-1}>x,\tau\ge n\}).
\end{eqnarray*}
Therefore,
\begin{eqnarray*}
{\bf P}\{S_\tau>x\} &=&
\overline F(x)
+\sum_{n=2}^N ({\bf P}\{S_{n-1}\le x, S_n>x,\tau\ge n\}
-{\bf P}\{S_{n-1}>x, S_n\le x,\tau\ge n\}).
\end{eqnarray*}
Now it suffices to show that, for each $n$,
\begin{eqnarray}\label{est.N.1}
P_1 \equiv {\bf P}\{S_{n-1}\le x, S_n>x,\tau\ge n\}
&\sim& \overline F(x){\bf P}\{\tau\ge n\}
\end{eqnarray}
and
\begin{eqnarray}\label{est.N.2}
P_2 \equiv {\bf P}\{S_{n-1}>x, S_n\le x,\tau\ge n\}
&=& o(\overline F(x)).
\end{eqnarray}

The subexponentiality of $F$ implies that,
for each $n\ge 2$,
\begin{eqnarray}\label{subexp.col.1}
{\bf P}\{S_n>x\} &\sim& n\overline F(x)
\quad\mbox{ as }x\to\infty.
\end{eqnarray}
In particular, there exists $c$ such that, for
all $n=2,\ldots , N$,
\begin{eqnarray}\label{subexp.col.2}
{\bf P}\{S_n>x\} &\le& c\overline F(x)
\quad\mbox{ for all }x.
\end{eqnarray}
The subexponentiality of $F$ also implies,
for any $A(x)\to\infty$ such that
$\overline F(x+A(x))\sim\overline F(x)$,
\begin{eqnarray}\label{subexp.col.3}
\int_{A(x)}^{x+A(x)}\overline F(x-y)F(dy)
&=& o(\overline F(x))
\quad\mbox{ as }x\to\infty.
\end{eqnarray}

To establish (\ref{est.N.1}),
we first note that
$\{\tau\ge n\}=\overline{\{\tau\le n-1\}}$
and thus $\sigma(S_{n-1},{\bf I}\{\tau\ge n\})$
does not depend on $\xi_n$,
since $\tau$ does not depend on the future.
This implies
\begin{eqnarray*}
P_1 &=& \int_0^\infty
{\bf P}\{S_{n-1}\in(x-y,x],
\xi_n\in dy,\tau\ge n\}\\
&=& \int_0^\infty
{\bf P}\{S_{n-1}\in(x-y,x],\tau\ge n\} F(dy).
\end{eqnarray*}
We use the following decomposition, $A>0$:
\begin{eqnarray}\label{P.1.i1-i3}
P_1 &=& \Biggl(\int_0^A+\int_A^{x+A}
+\int_{x+A}^\infty\Biggr)
{\bf P}\{S_{n-1}\in(x-y,x],\tau\ge n\} F(dy)
\nonumber\\
&\equiv& I_1+I_2+I_3.
\end{eqnarray}
By (\ref{subexp.col.1}) and by the
long-tailedness of $F$, for any fixed $A$,
\begin{eqnarray}\label{i1.est}
I_1 &\le& {\bf P}\{S_{n-1}\in(x-A,x]\}
=o(\overline F(x))\quad\mbox{ as }x\to\infty.
\end{eqnarray}
By (\ref{subexp.col.2}) and (\ref{subexp.col.3})
we get, for $A=A(x)\to\infty$,
\begin{eqnarray}\label{i2.est}
I_2 &\le& \int_A^{x+A}
{\bf P}\{S_{n-1}>x-y\} F(dy)\nonumber\\
&\le& c\int_A^{x+A}\overline F(x-y) F(dy)
= o(\overline F(x))\quad\mbox{ as }x\to\infty.
\end{eqnarray}
Uniformly in $y\ge x+A(x)$,
${\bf P}\{S_{n-1}\in(x-y,x],\tau\ge n\}
\to{\bf P}\{\tau\ge n\}$ as $x\to\infty$.
Thus,
\begin{eqnarray}\label{i3.est}
I_3 &\sim&
{\bf P}\{\tau\ge n\}\overline F(x+A(x))
\sim {\bf P}\{\tau\ge n\}\overline F(x)
\quad\mbox{ as }x\to\infty.
\end{eqnarray}
Substituting (\ref{i1.est})--(\ref{i3.est})
into (\ref{P.1.i1-i3}) we get (\ref{est.N.1}).

To prove (\ref{est.N.2}) we note that
\begin{eqnarray*}
P_2 &\le& {\bf P}\{S_{n-1}\in(x,x+A]\}
+{\bf P}\{S_{n-1}>x+A\}F(-A).
\end{eqnarray*}
As in (\ref{i1.est}), the first term on the
right is of order $o(\overline F(x))$.
Due to (\ref{subexp.col.2}), the second term
is not greater than $c\overline F(x)F(-A)$
where $F(-A)$ can be made as small as we please
by the choice of sufficiently large $A$.
The proof is complete.

Here is our general result for a counting random
variable with, possibly, unbounded support.

\begin{Theorem}\label{thm.stopping.}
Let ${\bf E}|\xi|<\infty$ and let a counting
variable $\tau$ do not depend on the future.
Assume that $F\in{\mathcal S}^*$ and that there
exists an increasing function $h(x)$ such that
\begin{eqnarray}\label{hx.cond}
\overline F(x\pm h(x))\sim\overline F(x)
&\mbox{and}& {\bf P}\{\tau>h(x)\}
=o(\overline F(x))\quad\mbox{ as }x\to\infty.
\end{eqnarray}
Then ${\bf P}\{S_\tau>x\}\sim{\bf E}\tau\overline F(x)$
as $x\to\infty$.
\end{Theorem}

\proof[Proof\ \ {\rm of Theorem \ref{thm.stopping.}
follows from Lemmas \ref{lm.stopping.+}
and \ref{lm.stopping.-} below}]
Condition (\ref{hx.cond}) is stronger than
condition (\ref{eq1}).
At the end of this section, we provide an
example of a stopping time which shows
that condition (\ref{hx.cond}) is essential
and cannot be weakened to (\ref{eq1}).

\begin{Lemma}\label{lm.stopping.+}
Let ${\bf E}\xi>0$ and let a counting variable
$\tau$ do not depend on the future.
If $F$ is long-tailed then
\begin{eqnarray*}
\liminf_{x\to\infty}
\frac{{\bf P}\{S_\tau>x\}}{\overline F(x)}
&\ge& {\bf E}\tau.
\end{eqnarray*}
If, in addition, $F\in{\mathcal S}^*$ and
condition {\rm(\ref{hx.cond})} holds, then
${\bf P}\{S_\tau>x\}\sim{\bf E}\tau\overline F(x)$
as $x\to\infty$.
\end{Lemma}

\proof
Fix a positive integer $N$ and a positive $A$.
The following lower bound holds, for $x>A$:
\begin{eqnarray*}
{\bf P}\{S_\tau>x\} &\ge& \sum_{j=1}^N
{\bf P}\{S_1,\ldots,S_{j-1}\in[-A,A], \xi_j>x+2A,
S_\tau>x,\tau\ge j\}\\
&\ge& \sum_{j=1}^N
{\bf P}\{S_1,\ldots,S_{j-1}\in[-A,A], \xi_j>x+2A,
\min_{i>j}(S_i-S_j)>-A,\tau\ge j\}.
\end{eqnarray*}
Since $\{\tau\ge j\}=\overline{\{\tau\le j-1\}}$ and
since $\tau$ does not depend on the future,
\begin{eqnarray*}
{\bf P}\{S_\tau>x\} &\ge& \sum_{j=1}^N
{\bf P}\{S_1,\ldots,S_{j-1}\in[-A,A],\tau\ge j\}
{\bf P}\{\xi_j>x+2A,\min_{i>j}(S_i-S_j)>-A\}\\
&=& \overline F(x+2A)
{\bf P}\{\min_{i\ge 1}S_i>-A\} \sum_{j=1}^N
{\bf P}\{S_1,\ldots,S_{j-1}\in[-A,A],\tau\ge j\}.
\end{eqnarray*}
By the long-tailedness of $F$,
\begin{eqnarray*}
\liminf_{x\to\infty}
\frac{{\bf P}\{S_\tau>x\}}{\overline F(x)}
&\ge& {\bf P}\{\min_{i\ge 1}S_i>-A\} \sum_{j=1}^N
{\bf P}\{S_1,\ldots,S_{j-1}\in[-A,A],\tau\ge j\}.
\end{eqnarray*}
Since the mean of $\xi$ is positive,
${\bf P}\{\min_{i\ge 1}S_i>-A\}\to 1$ as $A\to\infty$.
Hence, for any $N$,
\begin{eqnarray*}
\liminf_{x\to\infty}
\frac{{\bf P}\{S_\tau>x\}}{\overline F(x)}
&\ge& \sum_{j=1}^N {\bf P}\{\tau\ge j\}.
\end{eqnarray*}
Letting now $N\to\infty$ completes the proof
of the lower bound.

The upper bound,
\begin{eqnarray*}
\limsup_{x\to\infty}
\frac{{\bf P}\{S_\tau>x\}}{\overline F(x)}
&\le& {\bf E}\tau,
\end{eqnarray*}
follows from [\ref{FZ_2003}, Corollary 3] which states that,
under the conditions
$F\in{\mathcal S}^*$ and (\ref{hx.cond}),
${\bf P}\{M_\tau>x\} \sim
\overline F(x){\bf E}\tau$ as $x\to\infty$.
The proof is complete.

\begin{Lemma}\label{lm.stopping.-}
Let ${\bf E}\xi\le0$ and let a counting variable
$\tau$ do not depend on the future.
If $F\in{\mathcal S}^*$, then
\begin{eqnarray*}
\limsup_{x\to\infty}
\frac{{\bf P}\{S_\tau>x\}}{\overline F(x)}
&\le& {\bf E}\tau.
\end{eqnarray*}
Under the additional condition {\rm(\ref{hx.cond})},
${\bf P}\{S_\tau>x\}\sim{\bf E}\tau\overline F(x)$
as $x\to\infty$.
\end{Lemma}

\proof
The upper bound follows from
[\ref{FZ_2003}, Corollary 3] in the same way as
the upper bound in the previous proof.
To obtain the lower bound, take any positive
$\varepsilon$ and consider a random walk
$\widetilde S_n=S_n+n (|{\bf E}\xi|+\varepsilon )$ with
a positive drift. We have
\begin{eqnarray*}
{\bf P}\{S_\tau>x\} &=&
{\bf P}\{\widetilde S_\tau>x+(|{\bf E}\xi|+\varepsilon )\tau\}\\
&\ge& {\bf P}\{\widetilde S_\tau>x+(|{\bf E}\xi|
+\varepsilon )h(x)\}
-{\bf P}\{\tau>h(x)\}.
\end{eqnarray*}
Here the last term in the right side is $o(\overline F(x))$ and,
by Lemma \ref{lm.stopping.+},
the first term is equivalent to
${\bf E}\tau\overline F(x+(|{\bf E}\xi|+\varepsilon )h(x))
\sim{\bf E}\tau\overline F(x)$ as $x\to\infty$.
This completes the proof.

For intermediate regularly varying tail distributions,
Theorem \ref{thm.stopping.} implies the following

\begin{Corollary}\label{cor.stopping.reg}
Let ${\bf E}|\xi|<\infty$ and let
a counting variable $\tau$ do not depend
on the future. Assume that $F$ is an intermediate
regularly varying distribution and that
\begin{eqnarray}\label{tau=osmallF}
{\bf P}\{\tau>x\} &=& o(\overline F(x))
\quad\mbox{ as }x\to\infty.
\end{eqnarray}
Then
${\bf P}\{S_\tau>x\}\sim{\bf E}\tau\overline F(x)$
as $x\to\infty$.
\end{Corollary}

The latter corollary generalises the corresponding
result by Greenwood and Monroe
[\ref{GreenwoodMonroe}, Theorem 1] where
a regularly varying $F$ and a stopping time
$\tau$ were considered.
In [\ref{BU}, Theorem 2], A. Borovkov and Utev
obtained an upper bound for the tail distribution
of $S_\tau$ assuming that both tail distributions
of $\xi_1$ and of $\tau$ are bounded from above
by the same dominated varying distribution.

\proof
From condition (\ref{tau=osmallF}),
for any $\varepsilon>0$,
\begin{eqnarray*}
{\bf P}\{\tau>\varepsilon x\} &=&
o(\overline F(\varepsilon x))=o(\overline F(x))
\quad\mbox{ as }x\to\infty,
\end{eqnarray*}
since $F$ is intermediate regularly varying.
Thus, there exists an increasing function
$h(x)=o(x)$ such that
${\bf P}\{\tau>h(x)\}=o(\overline F(x))$
as $x\to\infty$. Again by the intermediate
regular variation of $F$, for any $h(x)=o(x)$,
$\overline F(x\pm h(x))\sim\overline F(x)$.
So, condition (\ref{hx.cond}) is fulfilled
and we can conclude the desired asymptotics
from Theorem \ref{thm.stopping.}.

We conclude with an example of a stopping time
$\tau$ showing that condition (\ref{hx.cond})
is essential for the conclusion of
Theorem \ref{thm.stopping.}.
Consider a distribution $F$ on $[1,\infty)$.
Take an increasing function
$H(x):{\bf R}\to{\bf Z}^+$ such that $H(x)<x/2$.
The counting random variable $\tau=H(2\xi_1)+1$
is a stopping time. On the event $\xi_1>x-H(x)$
we have $\tau\ge H(2(x-H(x)))+1\ge H(x)+1$. Hence,
\begin{eqnarray*}
{\bf P}\{S_\tau>x\} \ge
{\bf P}\{\xi_1>x-H(x),\xi_2+\ldots+\xi_\tau\ge H(x)\}
= {\bf P}\{\xi_1>x-H(x)\},
\end{eqnarray*}
due to $\xi\ge1$. For a Weibull type distribution,
namely $\overline F(x)=e^{-x^\beta}$,
$0<\beta<1$, $x\ge 1$,
we can choose $H(x)$ in such a way that $H(x)=o(x)$
and $H(x)/x^{1-\beta}\to\infty$ as $x\to\infty$.
Then condition (\ref{eq1}) holds,
but asymptotics (\ref{eq_lim}) does not,
because $\overline F(x-H(x))/\overline F(x)\to\infty$
and
\begin{eqnarray*}
\frac{{\bf P}\{S_\tau>x\}}{\overline F(x)}
&\to& \infty.
\end{eqnarray*}
In this example there is no a function $h(x)$ such that
condition (\ref{hx.cond}) holds. Indeed,
if $\overline F(x-h(x))\sim\overline F(x)$ then
$h(x)=o(x^{1-\beta})$ and $H^{-1}(h(x)-1)=o(x)$
which implies
\begin{eqnarray*}
{\bf P}\{\tau>h(x)\}/\overline F(x)
&=& {\bf P}\{H(2\xi)>h(x)-1\}/\overline F(x)
\to\infty \quad\mbox{ as }x\to\infty.
\end{eqnarray*}

\section*{\normalsize References}

\newcounter{bibcoun}
\begin{list}{\arabic{bibcoun}.}{\usecounter{bibcoun}\itemsep=0pt}
\small

\item\label{ALS}
Ale\v skevi\v cien\.e, A., Leipus, R. and \v{S}iaulys, J., 2008.
Tail behavior of random sums under consistent
variation with applications to the compound
renewal risk model.
Extremes 11, 261-–279.

\item\label{AN}
Athreya, K. B. and Ney, P. E., 1973.
Branching processes.
Springer.

\item\label{BB}
Borovkov, A. A. and Borovkov, K. A., 2008.
Asymptotic analysis of random walks:
Heavy-tailed distributions.
Cambridge University Press.

\item\label{BM}
Borovkov, A. A. and Mogul'skii, A. A., 2006.
Integro-local and integral theorems for sums of
random variables with semiexponential distributions.
Siberian Math. J. 47, 990--1026.

\item\label{BU}
Borovkov, A. A. and Utev, S. A., 1993.
Estimates for distributions of sums stopped at a Markov time.
Theory Probab. Appl. 38, 214-�225.

\item\label{DOV}
Daley, D. J., Omey, E. and Vesilo, R., 2007.
The tail behaviour of a random sum of subexponential
random variables and vectors.
Extremes 10, 21--39.

\item\label{DDS}
Denisov, D., Dieker, A. B. and Shneer, V., 2008.
Large deviations for random walks under
subexponentiality: the big-jump domain.
Ann. Probab., to appear.

\item\label{DFK2004}
Denisov, D., Foss, S. and Korshunov, D., 2004.
Tail asymptotics for the supremum of a random walk
when the mean is not finite.
Queueing Systems 46, 15--33.

\item\label{DFK1}
Denisov, D., Foss, S. and Korshunov, D., 2008.
On lower limits and equivalences for
distribution tails of randomly stopped sums.
Bernoulli 14, 391--404.

\item\label{DFK2}
Denisov, D., Foss, S. and Korshunov, D., 2007.
Lower limits for distribution tails of randomly stopped sums.
Theory Probab. Appl. 52, 794--802
(Russin; to be translated in 2008).


\item\label{EGV}
Embrechts, P., Goldie, C. M. and Veraverbeke, N., 1979.
Subexponentiality and infinite divisibility.
Z. Wahrscheinlichkeitstheorie verw. Gebiete
49, 335--347.

\item\label{FGMS}
Fa\"y, G., Gonz\'alez-Ar\'evalo, B.,
Mikosch, T. and Samordnitsky, G., 2006.
Modeling teletraffic arrivals by a Poisson cluster process"
Queueing Systems 54, 121--140.

\item\label{FK}
Foss, S. and Korshunov, D., 2007.
Lower limits and equivalences for convolution tails.
Ann. Probab. 35, 366--383.

\item\label{FPZ}
Foss, S., Palmowski, Z. and Zachary  S., 2005.
The probability of exceeding a high boundary on a random
time interval for a heavy-tailed random walk.
Ann. Appl. Probab. 15, 1936--1957.

\item\label{FZ_2003}
Foss, S. and Zachary, S., 2003.
The maximum on a random time interval of a random
walk with long-tailed increments and negative drift.
Ann. Appl. Probab. 13, 37--53.

\item\label{Greenwood}
Greenwood, P., 1973.
Asymptotics of randomly stopped sequences
with independent increments.
Ann. Probab. 1, 317--321.

\item\label{GreenwoodMonroe}
Greenwood, P. and Monroe, I., 1977.
Random stopping preserves regular variation
of process distributions.
Ann. Probab. 5, 42--51.

\item\label{Kl}
Kl\"uppelberg, C., 1988.
Subexponential distributions and integrated tails.
J. Appl. Probab. 25, 132--141.

\item\label{KP}
Kolmogorov, A. N. and Prokhorov, Yu. V., 1949.
On sums of a random number of random terms. (Russian).
Uspehi Matem. Nauk IV, no. 4, 168--172.

\item\label{K2002}
Korshunov, D., 2002.
Large-deviation probabilities for maxima of sums of
independent random variables with negative mean and
subexponential distribution.
Theory Probab. Appl. 46, 355--366.

\item\label{Linnik}
Linnik, Yu. V., 1961.
On probability of large deviations for the sums
of independent variables.
Proc. Fourth Berkeley Symp. Math. Statist.
Probability, Univ. Calif. Press, V. 2, 289--306.

\item\label{MAV}
Mikosch, T. and Nagaev, A. V., 1998.
Large deviations of heavy-tailed sums
with applications in insurance.
Extremes 1, 81--110.

\item\label{AN1}
Nagaev, A. V., 1969.
Limit theorems that take into account large
deviations when Cramer's condition is violated.
(Russian) Izv. Akad. Nauk UzSSR Ser. Fiz.-Mat. Nauk 13,
No. 6, 17--22.

\item\label{AN2}
Nagaev, A. V., 1969.
Integral limit theorems taking large deviations into
account when Cram\'er's condition does not hold I, II.
Theory Probab. Appl. 14, 51--64, 193--208.

\item\label{SN1960}
Nagaev, S. V., 1960.
Local limit theorems for large deviations. (Russian).
Theory Probab. Appl. 5, 259--261.


\item\label{SN1962}
Nagaev, S. V., 1962.
An integral limit theorem for large deviations. (Russian).
Izv. Akad. Nauk UzSSR Ser. Fiz.-Mat. Nauk, no. 6, 37--43.


\item\label{SN1973}
Nagaev, S. V., 1973.
Large deviations for sums of independent random variables.
Trans. Sixth Prague Conf. (1971).
Inform. Theory, Statist. Decision Functions,
Random Processes, 657--674.

\item\label{SN}
Nagaev, S. V., 1979.
Large deviations of sums of independent random variables.
Ann. Probab. 7, 745--789.

\item\label{SN1981}
Nagaev, S. V., 1981.
On asymptotic behaviour of probabilities of one-sided large deviations.
Theory Probab. Appl. 26, 362--366.

\item\label{Ng_etal}
Ng, K. W., Tang, Q. H., and Yang, H., 2002.
Maxima of sums of heavy-tailed random variables.
ASTIN Bull. 32, 43--55.

\item\label{Pitman}
Pitman, E. J. G., 1980.
Subexponential distribution functions.
J. Austral. Math. Soc. Ser. A 29,
337--347.

\item\label{Roz1990}
Rosovskii, L. V., 1990.
Probabilities of Large Deviations of Sums of Independent Random
Variables with Common Distribution Function in
the Domain of Attraction of the Normal Law.
Theory Probab. Appl. 34, 625--644.

\item\label{Roz1993}
Rosovskii, L. V., 1993.
Probabilities of large deviations on the whole axis.
Theory Probab. Appl. 38, 53--79.

\item\label{Roz1997}
Rosovskii, L. V., 1997.
Probabilities of large deviations for sums of independent
random variables with a common distribution function
from the domain of attraction of an asymmetric stable law.
Theory Probab. Appl. 42, 454--482.

\item\label{Rudin}
Rudin, W., 1973.
Limits of ratios of tails of measures.
Ann. Probab. 1, 982--994.

\item\label{Shneer}
Shneer, V., 2004.
Estimates for the distributions of the sums
of subexponential random variables.
Sib. Math. J. 45, 1143--1158.

\item\label{Stam}
Stam, A. J., 1973.
Regular variation of the tail of a
subordinated probability distribution.
Adv. Appl. Prob. 5, 308--327.

\item\label{Wachtel}
Wachtel, V., 2007.
Limit theorems for large deviation probabilities
for critical Galton-Watson process with power tails.
Theory Probab. Appl. 52, 644--659
(Russin; to be translated in 2008).

\end{list}

\noindent School of MACS,\\
Heriot-Watt University,\\
Edinburgh EH14 4AS, UK.\\
E-mail: Denisov@ma.hw.ac.uk\\
Foss@ma.hw.ac.uk\\[5mm]
Sobolev Institute of Mathematics\\
4 Koptyuga pr., Novosibirsk 630090, Russia.\\
E-mail: Korshunov@math.nsc.ru

\end{document}